\documentclass[12pt]{iopart}
\usepackage{iopams}

\newcommand{\RR}{\mathbb{R}}

\newcommand{\intp}[1]{
\int _{\RR} #1 \;\rmd p}

\newcommand{\intx}[1]{
\int _0^L #1 \;\rmd x}

\newcommand{\intxp}[1]{
\int _0 ^L \int _{\RR} #1 \;\rmd p\;\rmd x}

\newcommand{\intyq}[1]{
\int _0 ^L \int _{\RR} #1 \;\rmd q\;\rmd y}

\newtheorem{theorem}{Theorem}
\newtheorem{proposition}[theorem]{Proposition}
\newtheorem{remark}{Remark}

\newenvironment{proof}{\textit{Proof:\quad}}{\hfill$\square$}

\begin{document}

\title[Boltzmannian waves in laser-plasma interaction]{On the harmonic Boltzmannian waves in laser-plasma interaction}

\author{Mihai Bostan$^1$ and Simon Labrunie$^2$}

\address{$^1$ Laboratoire de Math\'ematiques de Besan{\c c}on, UMR CNRS 6623, Universit\'e de Franche-Comt\'e, 16 route de Gray, 25030 Besan{\c c}on  cedex,
France.}
\address{$^2$ Institut \'Elie Cartan (Math\'ematiques) UMR 7502, Nancy-Universit\'e, CNRS et INRIA (projet CALVI), 54056 Vand{\oe}uvre-l\`es-Nancy cedex, France.}
\ead{mbostan@math.univ-fcomte.fr, labrunie@iecn.u-nancy.fr}

\begin{abstract}
We study the permanent regimes of the reduced Vlasov--Maxwell
system for laser-plasma interaction. A non-relativistic and two different relativistic models are investigated. We prove the existence of solutions where the distribution function is Boltzmannian and the electromagnetic variables are time-harmonic and circularly polarized. 
\end{abstract}


%02.30.Hq Ordinary differential equations
%02.30.Ks Delay and functional equations
%02.30.Xx Calculus of variations
%52.27.-h Basic studies of specific kinds of plasmas
%52.27.Aj Single-component, electron-positive-ion plasmas
%52.35.Mw Nonlinear phenomena: waves, wave propagation, and other interactions (including parametric effects, mode coupling, ponderomotive effects, etc.)
%52.38.-r Laser-plasma interactions (for plasma production and heating by laser beams, see 52.50.Jm)

\submitto{\JPA}

%%%%%%%%%%%%
%          %
%   TEXT   %
%          %
%%%%%%%%%%%%

\section{Introduction.\label{intro}}
The in-depth understanding of laser-plasma interaction is of paramount importance for the eventual success of inertial confinement fusion research,
but is also interesting for magnetic confinement fusion research,
since tokamak plasmas can be heated by electromagnetic waves. 
The complex kinetic phenomena involved in this interaction, and the
instabilities they may generate~\cite{Alma05}
need to be studied by kinetic models~\cite{JBG+92}, even though
hydrodynamic models~\cite{CoCo04} are more affordable to simulate complex, 
high-dimensional geometries.
However, the use of the full 3D Vlasov--Maxwell system is of course
impossible in most practical situations. Therefore, 
the \emph{reduced Vlasov--Maxwell system for 
laser-plasma interaction} (hereafter called the ``laser-plasma system''; see~\mbox{(\ref{vla.adim}--\ref{ond.adim})} below)
was introduced in~\cite{JBG+92}. The model has been shown to capture some essential features of this interaction~\cite{JBG+92,Alma05}; and it has been successfully used for deriving relevant physical models in novel situations~\cite{LBK+02}.

\medbreak

The laser-plasma system has been the object of 
several mathematical investigations~\cite{CaLa06,Bost06,BosIUMJ06}.
In this framework, it is interesting to find classes of exact solutions
which may serve as ``reference solutions'', to which other solutions
may be compared in order to study the dynamic of the interaction.
Reference solutions for the Vlasov--Poisson system are, for instance, 
the Bernstein--Greene--Kruskal or \emph{BGK modes}~\cite{BeGK57}, 
given by a distribution function of the form~$f(W)$, where $W$~is 
the energy of one particle. When the function~$f$ is convex, as in the 
Maxwellian case $f(W) \propto \rme^{-W/\theta}$, such solutions 
represent fundamental equilibrium states; their existence and stability 
are well known~\cite{CaCD02}.
For a general~$f$, BGK~solutions may represent various wave phenomena; they have been the object of many investigations in the physical community~\cite[and references therein]{MaBe00}. The mathematical theory is still less developed; 
interesting existence and (in)stability results have appeared recently~\cite{MaBe00,BGK-recent}.
For both Vlasov--Poisson and Vlasov--Maxwell systems, there are also 
linearised solutions leading to the dispersion relations of the various
types of waves; \emph{e.g.} for electromagnetic waves 
$\omega^2 = \omega_\mathrm{p}^2 + k^2$, where
$(\omega,k)$ are the pulsation and wave number, and 
$\omega_\mathrm{p}$~is the plasma pulsation. 

\medbreak

In this article, we shall introduce a class of exact solutions to
the laser-plasma system which generalises, at the same time, 
Maxwellian equilibria and linear electromagnetic waves.
Indeed, we investigate the existence of
quasi-static solutions where the distribution function is at any
time proportional to the Boltzmann factor; this static character can be
reconciled with the electromagnetic character of the system by assuming
a harmonic time dependence of the electromagnetic field and a circular 
polarisation. This ansatz was already used in~\cite{LBK+02}, but in a
different physical and mathematical context. The latter work 
investigates the existence of solitons in an electron-positron plasma, 
where no charge separation occurs. Here we are dealing with a general 
ion-electron plasma, and we are looking for space periodic solutions.

\medbreak

The paper is organized as follows.  We recall the mathematical results known about the laser-plasma system and introduce the quasi-static model in section~\ref{deriv}. Then, in section~\ref{boltz} we solve (in the space periodic setting) the so-called \emph{Boltzmann problem}, which consists in finding
the equilibrium density given the electromagnetic potentials, and we estimate its solutions. 
In section~\ref{fixed} we construct a fixed point application
and we study its properties. The existence of Boltzmaniann
equilibria then follows by applying the Schauder fixed point theorem.
Several extensions of the model are briefly discussed in section~\ref{exten}, and we conclude in section~\ref{concl}.

\section{The harmonic Boltzmannian model.\label{deriv}}
%\subsection{Context.}
The reduced Vlasov--Maxwell system for laser-plasma interaction
describes the evolution of the distribution function of a population of 
electrons in a one space dimensional plasma interacting with a laser 
wave. In a first approach, we assume that the ions are at rest and their density is given --- which is physically acceptable at the time scale of a laser wave.
After a suitable rescaling~\cite{CaLa06}, this system can be cast 
in the following form:
\begin{eqnarray}
\frac{\partial f}{\partial t} +  \frac{p}{\gamma _1}\,
\frac{\partial f}{\partial x} - \left (E(t,x) +
\frac{\bi{A}(t,x)}{\gamma _2} \cdot \frac{\partial \bi{A}}{\partial x} 
\right)\, \frac{\partial f}{\partial p} = 0,
\label{vla.adim}\\
\frac{\partial E}{\partial x} = \rho_\mathrm{b} (x) - \rho
(t,x),\qquad \frac{\partial E}{\partial t } - j(t,x) = 0,
\label{poi.adim}\\
\frac{\partial^2 \bi{A}}{\partial t^2} - \frac{\partial^2
\bi{A}}{\partial x^2} + \tilde\rho(t,x)\,\bi{A} (t,x) = 0,
\label{ond.adim}
\end{eqnarray}
where: $f(t,x,p)$ is the electron distribution function
($p$~denotes the $x$-compo\-nent of the momentum vector);
$E$~is the $x$-component if the electric field;  $\bi{A} =
(0,A_y,A_z)$ is the vector potential of the laser wave; $\rho_\mathrm{b}(x)$ is the (static) background ion density; $\gamma _1, \gamma
_2 $ are Lorentz
factors. We distinguish three cases:
\begin{enumerate}
\item the non-relativistic case (NR), $\gamma _1 = \gamma _2 =1$;
\item the quasi-relativistic case (QR), $\gamma _1 = (1 + p^2)^{1/2}, \gamma _2 =1$;
\item the fully relativistic case (FR), $\gamma _1 = \gamma _2 = ( 1
+ p^2 + |\bi{A}|^2)^{1/2}$, which is the original model
of~\cite{JBG+92}.
\end{enumerate}
The moments $\rho,\ \tilde\rho,\ j$ are given by
\begin{equation}
\fl
\label{Equ5} \rho(t,x):= \intp{f(t,x,p)}, \quad \tilde\rho(t,x):=
\intp{\frac{f(t,x,p)}{\gamma _2}},\quad j(t,x):=
\intp{\frac{p}{\gamma _1}f(t,x,p)}.
\end{equation}
We supplement the system (\ref{vla.adim}, \ref{poi.adim}, 
\ref{ond.adim}) with initial conditions
\begin{equation}
\fl
\label{Equ6} f(0,x,p) = f_0 (x,p), \;(x,p) \in \RR ^2, \;\;(E, \bi{A},
\partial _t \bi{A}) (0,x) = (E_0, \bi{A}_0, \bi{A} _1),\;x\in
\RR.
\end{equation}
In~\cite{CaLa06} it was proved that, for suitable initial
conditions, (\ref{vla.adim}--\ref{Equ6}) has a
unique classical solution, which is global in time in the QR~case,
and local in time in the NR~case. In the latter case, the classical
solution can be extended to a global weak solution with
$f$~continuous and~$\bi{A}$ continuously differentiable in all their
variables. The FR~model was studied in~\cite{Bost06}. It was shown
that (\ref{vla.adim}--\ref{Equ6}) admits a unique global
classical solution preserving the total energy. The stationary
solutions of these models in a bounded domain have been analysed
in~\cite{BosIUMJ06}. 

\medbreak

All three models admit space periodic solutions. If the initial
data are \mbox{$L$-periodic} in~$x$ and satisfy the neutrality condition
\begin{equation}
\label{eq:neutr} \intxp{f_0(x,p)} = \intx{\rho_\mathrm{b} (x)}=:M,
\end{equation}
then, by using the continuity equation $\partial _t \rho + \partial _x j = 0$,
we deduce that the system remains globally neutral at any time $t
>0$
\begin{equation}
\label{Equ8bis} \intxp{f(t,x,p)} = \intxp{f_0(x,p)} =\intx{\rho_\mathrm{b}
(x)}.
\end{equation}
By uniqueness of the solution one gets also that $(f(t),E(t),\bi{A}(t))$ are \mbox{$L$-periodic} in space for any $t >0$. {F}rom now on, we work in
the framework of periodic functions: all differential equations will
be implicitly supplemented with \mbox{$L$-periodic} boundary conditions.

\smallbreak

{F}rom~\eref{Equ8bis} we deduce the existence of a unique function 
$V = V(t,x)$, satisfying 
$\partial _x ^2 V(t,x) = \rho_\mathrm{b} (x) - \rho (t,x)$, $V(t,0) = 0$ and
$(V,\partial _x V)(t,x) = (V,\partial _x V)(t, x +L)$, for all $(t,x)
\in [0,+\infty )\times \RR $.
The field~$E$ derives from the potential~$V$, {\it i.e.}, $E = \partial _x V$. 

\medbreak

The purpose of this article is to study the existence of
particular solutions of~(\ref{vla.adim}, \ref{poi.adim}, 
\ref{ond.adim}) corresponding to local Boltzmannian equilibria. 
These are defined by $f(t,x,p) \propto \rme^{-W(t,x,p)/\theta}$
where $W(t,x,p)$ is the energy of one particle being at the phase 
space point~$(x,p)$ at time~$t$, and $\theta $ is the scaled temperature. 
As it is well known, such functions are solutions to the Vlasov 
equation~\eref{vla.adim} iff $W$~is independent of time.
Thus, we assume that $V$ does not depend on~$t$, and that $\bi{A}$ is
time-harmonic and circularly polarized, {\it i.e.}, 
$$A_y (t,x) + \rmi\, A_z (t,x) = a(x) \, \rme^{\rmi\omega t},\quad 
\mbox{with \textit{a~priori}}\quad a(x) \in \mathbb{C}.$$
Then the energy $W(x,p)$ is given, according to the relativistic 
character, by
\begin{eqnarray*}
W(x,p) = \case12\, (p^2 + |a(x)|^2) + V(x),\;\;\mbox{in the
NR case},\\
W(x,p) = \sqrt{1+p^2} + \case12\, |a(x)|^2 + V(x),\;\;\mbox{in the
QR case},\\
W(x,p) = \sqrt{1+p^2 + |a(x)|^2}  + V(x),\;\;\mbox{in the FR case}.
\end{eqnarray*}
Imposing the constraint
$(\ref{eq:neutr})$ yields
\begin{equation}
\label{Equ9}
 f(x,p) = M\,
\frac{\rme^{-W(x,p)/\theta}}{\intyq{\rme^{-W(y,q)/\theta}}},
\quad \forall(x,p)\in \RR ^2.
\end{equation}
By direct computation we check that in all three cases $f$~solves
the Vlasov equation~\eref{vla.adim}. We then observe that 
$j(x) = \intp{\frac{p}{\gamma _1}f(x,p)}= 0$, for $x\in \RR$, 
and thus the system~(\ref{vla.adim}, \ref{poi.adim}, \ref{ond.adim}) 
reduces to
\begin{eqnarray}
\label{Equ10} V''(x) = \rho_\mathrm{b} (x) - \rho (x),&& x \in \RR, \\
\label{Equ11} - \omega ^2 a(x) - a''(x) = - \tilde\rho (x) a(x),&\quad& x
\in \RR.
\end{eqnarray}
with $\rho = \intp{f}$, $\tilde\rho = \intp{\frac{f}{\gamma _2}}$
and $f$ given by~\eref{Equ9}.

\medbreak

Of course, we are interested in solutions such that $a \not\equiv0$,
otherwise we find a Vlasov--Poisson equilibrium. If such a solution exists,
$a$~appears as an eigenfunction of the operator $A_{\tilde\rho} :=
-\frac{\rmd^2}{\rmd x^2} + \tilde\rho(x)$, associated to the eigenvalue~$\omega^2$. 
It is well known that these eigenvalues are real and generically simple;
in particular, the lowest eigenvalue is always simple. As the 
coefficients of~$A_{\tilde\rho}$ are real, we infer that both $\Re(a)$ and 
$\Im(a)$ are eigenfunctions; thus, generically, they must be 
proportional. In other words, $a(x) = \mathsf{a}(x)\, \rme^{\rmi \varphi}$,
where $\mathsf{a}$ is a \emph{real} eigenfunction and $\varphi\in\RR$.
Then, $|a(x)|^2 = \mathsf{a}(x)^2$, and $W,\ f,\ \rho,\ \tilde\rho$ are 
defined in terms of~$\mathsf{a}$; while we may take $\varphi=0$ by 
rotating the axes $Oy,~Oz$. This means that, without loss of generality,
we may restrict our search to \emph{real} functions~$a$ solution 
to~\eref{Equ11}.

\medbreak

We now rewrite the model~(\ref{Equ10}, \ref{Equ11}) in a form which will prove more convenient for analysis. We shall denote by the subscript~$\#$ the spaces of \mbox{$L$-periodic} functions, \emph{e.g.}:
$L^1_{\#}(\RR):= \{g \in L^1_{\mathrm{loc}}(\RR) : \forall x,\; g(x+L) = g(x)\},
\ C^0_{\#}(\RR){} := \{w \in C^0 (\RR): \forall x,\; w(x+L) = w(x)\}$.
First, we introduce the operator $\Phi : L^1_{\#}(\RR) \rightarrow C^0_{\#}(\RR)$ given by
\begin{equation*}
\Phi [g] = w \in C^0_{\#}(\RR){},\;\; - w''(x) = g(x),\;\;x\in (0,L),
\;\;w(0)= w(L) = 0,
\end{equation*}
for any $g \in L^1_{\#}(\RR)$. 
Then, we consider the function $\psi : \RR \rightarrow \RR$ given by
\begin{equation*}
\rme^{-\psi (x)/\theta} = \intp{\exp{- \frac{W(x,p) -
V(x)}{\theta}}}, \quad x \in \RR,
\end{equation*}
namely, according to the relativistic character:
\begin{eqnarray}
\mbox{NR:}\quad \rme^{-\psi (x)/\theta} 
&=& \rme^{- a(x)^2 / 2\theta }
\intp{\rme^{- p^2 /2 \theta }},
\label{df:psi:NR}\\
\mbox{QR:}\quad \rme^{-\psi (x)/\theta} 
&=& \rme^{- a(x)^2 / 2\theta }
\intp{\exp{- \frac{\sqrt{1 + p^2}}{ \theta}}},
\label{df:psi:QR}\\
\mbox{FR:}\quad \rme^{-\psi (x)/\theta} 
&=& \intp{\exp{- \frac{\sqrt{1 + p^2 +
a(x)^2}}{ \theta}}}.
\label{df:psi:FR}
\end{eqnarray}
Notice that there is a constant $C(\theta) \in \RR$ such that
\begin{equation}
\label{EquNRQR} \psi (x)= \frac{a(x)^2}{2} + C(\theta),\;\;\mbox{in the NR
and QR cases}.
\end{equation}
In the FR case, by observing that
$$
\case{1}{2}( \sqrt{1 + p^2} + |a(x)|) \leq \sqrt{ 1 + p^2 +
a(x)^2} \leq \sqrt{1 + p^2} + |a(x)|,
$$
we obtain
$$
C_2(\theta) \;\rme^{- |a(x)| / \theta} \leq \rme^{- \psi(x) / \theta} \leq C _1(\theta) \;\rme^{- |a(x)| / 2\theta },
$$
with
$$
C _1(\theta) := \intp{\exp {- \frac{\sqrt{1 + p^2}}{2\theta}}} >
\intp{\exp {- \frac{\sqrt{1 + p^2}}{\theta}}}=:C_2(\theta).
$$
Finally one gets
\begin{equation}
\label{EquFR} \frac{|a(x)|}{2} \leq \psi (x) + \theta \ln
C_1(\theta) \leq |a(x)| + \theta  \ln
\frac{C_1(\theta)}{C_2(\theta)}.
\end{equation}

\medbreak

The density $\rho$ can be expressed in function of~$\psi$; and the system~(\ref{Equ9}--\ref{Equ11}) can be recast as:
\begin{eqnarray}
f(x,p) = K\, \rme^{-W(x,p)/\theta}, \nonumber\\
\label{Equ13} \rho (x) = K \rme^{- \frac{\psi (x) +
V(x)}{\theta}},\;\;x \in \RR, \;\;V = \Phi [\rho - \rho_\mathrm{b}], \\
\label{Equ13bis} \tilde\rho(x) = K\, \intp{\frac{1}{\gamma _2}{\rme^{-W(x,p)/\theta}}}, \\
\label{Equ14} - \omega ^2 a(x) - a''(x) = -\tilde\rho(x)\,a(x),
\end{eqnarray}
where the constant $K = M \left ( \intyq{\rme^{- W(y,q)/\theta}}
\right)^{-1}$ is such that $\intx{\rho (x)} = M$. We call this
system the \emph{Boltzmann--Helmholtz equations}; they
can be seen as a sort of non-linear eigenvalue problem.

\section{The Boltzmann problem.\label{boltz}}
For the moment we suppose that the function $\psi $ is given and
we solve the so-called \emph{Boltzmann problem}~\eref{Equ13}. 
The proof of the following proposition is immediate and left to
the reader.
\begin{proposition}
\label{EstimPhi} For any function $g \in L^1_{\#}(\RR)$ we have:
$$
\|\Phi [g]\|_{L^\infty(\RR){}} \leq L \|g \|_{L^1(0,L){}}.
$$
If the function $g$ satisfies $\int _0 ^L g(x)\;\rmd x = 0$, then
$\Phi [g] \in C^1_{\#}(\RR)$ and we have:
$$
\left \| \frac{\rmd}{\rmd x}\Phi [g]\right \| _{L^\infty(\RR){}} \leq \|g
\|_{L^1(0,L){}}.
$$
\end{proposition}

\begin{proposition}
\label{Boltzmann} Let $\psi \in L^\infty_{\#}(\RR){}$, $u_\mathrm{b} \in L^1_{\#}(\RR),\ u_{\mathrm{b}}
\geq 0$, $M = \intx{u_{\mathrm{b}} (x)}$ and $\theta >0$. Then there is a
unique function $u \in L^1_{\#}(\RR)$ such that:
\begin{equation}
\label{Equ15} u = M\, \frac{\exp {- \frac{\psi + \Phi
[u-u_{\mathrm{b}}]}{\theta}}}{\int_0 ^L \exp {- \frac{\psi (y) + \Phi [u -
u_{\mathrm{b}}](y)}{\theta}}\;\rmd y}.
\end{equation}
Moreover it satisfies:
\begin{equation}
\fl\mbox{}\qquad
0 \leq u \leq  \inf_{C\in\RR} \frac{M}{L} \exp {\frac{1}{\theta}
\left ( \frac{1}{L}\int _0 ^L (\psi (y)-C)\, \rmd y - \inf_{\RR}
(\psi-C)  + 4\,L\,M\right )} =: u _{\psi} \,;
\label{Equ16a}
\end{equation}
and if $\psi \in W ^{1,\infty }(\RR)$ then $u \in W ^{1,\infty }(\RR)$ and:
\begin{equation}
\label{Equ16b} \mathrm{Lip} \;u \leq \frac{\mathrm{Lip} \;\psi +
2\,M}{\theta } u _{\psi}.
\end{equation}
\end{proposition}
\begin{proof}
One readily checks that~\eref{Equ15} is equivalent to the
minimization of the functional
$$
J[v] := \intx{ \{ \theta \sigma (v(x)) + \frac{1}{2}\left|
\frac{\rmd}{\rmd x}\Phi [v - u_{\mathrm{b}}]\right |^2 + \psi (x)v(x) \}},
$$
under the constraint $\intx{v(x)} = M$, where $\sigma (s) = s\ln
s$, $s >0$ and $\sigma (0) = 0$. This problem is a variant of that
considered in~\cite{CaCD02,Dol99} and its well-posedness
follows from a similar argument. The functional~$J$ is strictly
convex, l.s.c. and bounded from below on the set
$$
\mathcal{K}(L,M) = \left\{ v\in L^1_{\#}(\RR)\; : \; v \geq 0,\;\;
\intx{v(x)} = M \right\}.
$$
Indeed, by applying the Jensen inequality with the convex function
$\sigma $, the measure $\rmd\mu = \rme^{- \psi (x)/\theta}
\left ( \int _0 ^L \rme^{- \psi (y) /\theta }\;\rmd y \right
)^{-1} \rmd x $ and the function $v / \rme^{-\psi/\theta}$, one gets:
$$
J[v] \geq \intx{ \{\theta \sigma (v(x)) + \psi (x) v(x)\}} \geq
\theta M \ln \left[ M\, \left(\int _0 ^L \rme^{-\psi(y)/\theta}\;\rmd y \right)^{-1} \right],
$$
saying that $\inf _{v \in \mathcal{K}(L,M)} J[v] > -\infty $. Take
a minimising sequence $(u_n)_n$. By using the Dunford--Pettis
criterion we can assume (after a suitable extraction) that
$(u_n)_n$ converges weakly in $L^1(0,L)$ towards a function $u \in
\mathcal{K}(L,M)$. Since $J$ is convex we can pass to the limit by
involving the semi-continuity of $J$ and we obtain that $J[u]
=\inf _{v \in \mathcal{K}(L,M)} J[v]$. Writing the Euler--Lagrange
equation we obtain
$$
\theta ( 1 + \ln u ) + \Phi [u-u_{\mathrm{b}}] + \psi - \alpha = 0,
$$
where $\alpha $ enters as the Lagrange multiplier associated to
the constraint $\intx{u(x)} = M$ and thus we deduce~\eref{Equ15}. By using now the Jensen inequality with the
convex function $t \mapsto \rme^{-t}$, the measure $\rmd\mu =
L^{-1}\,{\rmd x}$ and the function $\left(\psi + \Phi[u-u_{\mathrm{b}}] \right)/ \theta$ we obtain:
$$
\exp\left[ -\frac{1}{L} \int _0 ^L \frac{\psi + \Phi
[u-u_{\mathrm{b}}]}{\theta}\; \rmd y\right] \leq \int _0 ^L \exp \left(- \frac{\psi  + \Phi
[u-u_{\mathrm{b}}]}{\theta}\right) \;\frac{\rmd y}{L}.
$$
Therefore  by using Proposition~\ref{EstimPhi} we deduce
\begin{equation*}
\fl
\left( \int _0 ^L \rme^ {- \frac{\psi + \Phi [u - u_\mathrm{b}]}{\theta}}\rmd y \right)^{-1} \leq \frac{1}{L}\, \exp {\frac{1}{L}\int _0 ^L
\rme^ {- \frac{\psi + \Phi [u - u_\mathrm{b}]}{\theta}}\rmd y} \leq \frac{1}{L}\,
\exp{\frac{1}{\theta} ( \frac{1}{L}\int _0 ^L \psi \;\rmd y +
2LM)},
\end{equation*}
which implies
\begin{equation*}
0 \leq u (x)  \leq   \frac{M}{L} \exp {\frac{1}{\theta}\, \left(
\frac{1}{L}\int _0 ^L \psi(y)\, \rmd y - \psi (x) + 4\,L\,M \right)}.  
\end{equation*}
This inequality is unchanged by replacing $\psi$ with $\psi-C$, for
any~$C\in\RR$; one thus infers~\eref{Equ16a}.
Assume now that $\psi \in W^{1, \infty }(\RR)$. By taking the
derivative with respect to $x$ in $(\ref{Equ15})$ one gets by
using Proposition~\ref{EstimPhi}
\begin{equation*}
|u'(x)|  = |u(x)|\;\left |\frac{\psi '(x) + \frac{\rmd}{\rmd x}\Phi
[u-u_{\mathrm{b}}]}{\theta}\right | \leq \|u \|_{L^\infty(\RR){}} \frac{\|\psi '
\|_{L^\infty(\RR){}} + 2\,M}{\theta},
\end{equation*}
and $(\ref{Equ16b})$ follows immediately.
\end{proof}

\section{The fixed point application.\label{fixed}}
For any $\mathfrak{a} > 0$ we define the fixed point application
${\cal F} _{\mathfrak{a}} : W^{1,\infty}_{\#} ( \RR) \rightarrow
W^{1,\infty}_{\#} ( \RR)$, ${\cal F}_{\mathfrak{a}}a = \tilde{a}$
for any $a \in
W^{1,\infty}_{\#} ( \RR)$ where:
\begin{itemize}
\item $\psi$ is given, according to the case, 
by~\eref{df:psi:NR}, \eref{df:psi:QR} or~\eref{df:psi:FR}; 
\item $\rho$ is the unique solution to the Boltzmann problem
$$
\rho = K e ^{- \frac{\psi + \Phi [\rho - \rho_\mathrm{b} ]}{\theta}},\;\;
\intx{\rho (x)} = M~;
$$
\item $ \tilde\rho = \rho $ in the NR and QR cases, while in the FR case:
$$
\tilde\rho (x) = \frac{M}{\int _0 ^L \exp {- \frac{\psi (y)+ \Phi [\rho
- \rho_\mathrm{b}](y)}{\theta}}\rmd y}\intp{\frac{ \exp {- \frac{\sqrt{1 + p^2 +
a(x)^2} + \Phi [\rho - \rho_\mathrm{b}](x)}{\theta}}}{\sqrt{1 + p^2+
a(x)^2}}}~;
$$
\item $\lambda $ is the first eigenvalue of the operator $A_{\tilde\rho} = -
\frac{\rmd^2}{\rmd x^2} + \tilde\rho $ with $L$ periodic boundary conditions,
\emph{i.e.},
\begin{equation}
\lambda = \inf _{b \in H^1_{\#}(\RR), b\neq 0} \frac{\intx{b'(x)^2 +
\tilde\rho (x)\, b(x)^2}}{\intx{b(x)^2}}~;
\label{Equ:lambda}
\end{equation}
\item $\tilde{a} $ is the corresponding eigenfunction of~$A_{\tilde\rho}$:
\begin{equation}
\label{Equ20} - \tilde{a}''(x) + \tilde\rho (x)\, \tilde a(x) = \lambda\,
\tilde{a} (x),\;\;x \in (0,L),
\end{equation}
\begin{equation}
\label{Equ21} \tilde{a}(0) = \tilde{a}(L),\;\;\tilde{a}'(0) =
\tilde{a}'(L),
\end{equation}
such that $\tilde{a} > 0$ and $\intx{\tilde{a}(x)^2} =
\mathfrak{a}^2$.
\end{itemize}
\begin{remark}
It is well known that the first eigenvalue of $A_{\tilde\rho}$ with $L$-periodic boundary conditions is simple~\cite{CodLev55} and that the eigenfunction  vanishes nowhere. Therefore $\tilde{a} = {\cal F}_{\mathfrak{a}}a$ is well defined.
\end{remark}

\medbreak

The properties of the application ${\cal F}_{\mathfrak{a}}$ are
summarized up below.
\begin{proposition}
\label{FixedPointAppli} Assume that $\rho_\mathrm{b} \in L^1_{\#}(\RR)$, $\rho_\mathrm{b}
\geq 0$, $\intx{\rho_\mathrm{b}(x)} = M$ and let $\mathfrak{a}, \theta$ be
positive real numbers. For any $a \in W^{1,\infty}_{\#}(\RR)$ such that
$\intx{a(x)^2} \leq \mathfrak{a}^2$ construct $\psi,\ \rho,\ \tilde\rho,\ 
\lambda$ and $\tilde{a} = {\cal F}_{\mathfrak{a}}a$ as above.
\begin{enumerate}
\item There are  constants $\rho _\star, a_\star$ depending on
$\mathfrak{a}, L, M, \theta $ such that
\begin{eqnarray*}
\fl \|\tilde\rho \|_{L^\infty(\RR){}} \leq \|\rho \|_{L^\infty(\RR){}} \leq \rho
_\star,\qquad \|\rho '\|_{L^\infty(\RR){}} \leq \frac{\|a\|_{L^\infty(\RR){}}
\|a'\|_{L^\infty(\RR){}} + 2\,M}{\theta}\rho _\star, \\
\fl \|\tilde\rho '\|_{L^\infty(\RR){}} \leq \frac{\|a\|_{L^\infty(\RR){}}
\|a'\|_{L^\infty(\RR){}} ( 1 + \theta) + 2\,M}{\theta} \rho _\star,\;\;0
\leq \lambda \leq \rho _\star,\;\;\|\tilde{a}\|_{W^{2,\infty}(\RR)}\leq
a_\star.
\end{eqnarray*}
\item ${\cal F}_{\mathfrak{a}}$ is continuous with respect to the
topology of $C^0_{\#}(\RR){}$ on the set ${\cal C} = \{a \in
C^0_{\#}(\RR){}\;:\; \|a\|_{L^2(0,L)} \le \mathfrak{a} \mbox{ and }
\|a\|_{L^\infty(\RR)} + \|a'\|_{L^\infty(\RR)} \leq a_\star \}$.
\end{enumerate}
\end{proposition}
\begin{proof}
(i) Take $a \in W^{1,\infty}_{\#}(\RR)$ such that $\|a\|_{L^2(0,L){}} \leq
\mathfrak{a}$. In the NR and QR cases, we deduce from~\eref{EquNRQR}
and~\eref{Equ16a} the bound
$$
0 \leq \rho \leq \frac{M}{L}\, \exp{\frac{1}{\theta}\left (
\frac{1}{L}\intx{\frac{a(x)^2}{2}} + 4\,L\,M \right )} \leq
\frac{M}{L}\, \exp {\frac{1}{\theta}\left (
\frac{1}{2L}\mathfrak{a}^2 + 4\,L\,M \right )}.
$$
In the FR case, combining~\eref{EquFR} and~\eref{Equ16a}
yields
\begin{eqnarray}
\fl 0 \leq \rho & \leq & \frac{M}{L}\, \exp {\frac{1}{\theta} \left (
\frac{1}{L} \int _0 ^L(\psi (y) + \theta \ln C_1 (\theta))\;\rmd y -
\inf ( \psi +\theta \ln C_1 (\theta)) + 4\,L\,M\right )} \nonumber \\
\fl & \leq & \frac{M}{L}\, \exp  {\frac{1}{\theta} \left( \frac{1}{L} \int
_0 ^L \left( |a(x)| + \theta \ln \frac{C_1(\theta)}{C_2(\theta)}\right)
\;\rmd y + 4\,L\,M\right)}\nonumber \\
\fl & \leq & \frac{M}{L}\, \exp  {\frac{1}{\theta} \left (
\frac{1}{\sqrt{L}} \mathfrak{a} +  \theta \ln
\frac{C_1(\theta)}{C_2(\theta)} + 4\,L\,M\right )} .\nonumber
\end{eqnarray}
We check easily that in all three cases we have $|\psi '(x)| \leq
|a(x)|\;|a'(x)|$, $x \in \RR$ and thus, by Proposition~\ref{Boltzmann} we deduce
$$\|\rho ' \|_{L^\infty(\RR){}} \leq \frac{\|a\|_{L^\infty(\RR){}}
\|a'\|_{L^\infty(\RR){}} + 2\,M}{\theta}\, \|\rho \|_{L^\infty(\RR){}}.$$
The estimate for~$\tilde\rho$ follows, since in all three cases 
$\gamma_2\ge1$ and $0\le\tilde\rho(x)\le\rho(x)$.
By taking the derivative with respect to $x$ in the expression of
$\tilde\rho$ we obtain by direct computation
$$
\|\tilde\rho ' \|_{L^\infty(\RR){}} \leq \left ( \|a\|_{L^\infty(\RR){}}
\;\|a'\|_{L^\infty(\RR){}} \left ( 1 + \frac{1}{\theta}\right ) +
\frac{2\,M}{\theta}\right ) \|\rho \|_{L^\infty(\RR){}}.
$$
We now estimate the eigenvalue~$\lambda$ and the
eigenfunction $\tilde{a}$. \Eref{Equ:lambda}~shows that $\lambda\ge0$ and,
by taking $b=1$, $\lambda\le\rho _\star$. Then, from~\eref{Equ20} we deduce
$$
\intx{\{\tilde{a}'(x)^2 + \tilde\rho (x)\, \tilde{a}(x)^2 \}} = \lambda
\intx{\tilde{a}(x)^2},
$$
and hence $\|\tilde{a}\|_{H^1(0,L)} \leq \sqrt{\lambda +1} \;\mathfrak{a}$. 
By using the Sobolev inclusion
$H^1(0,L) \subset L^\infty(0,L)$ one gets easily that
$\|\tilde{a}\|_{L^\infty(\RR)} + \|\tilde{a}'\|_{L^\infty(\RR)}  +
\|\tilde{a}''\|_{L^\infty(\RR)} \leq a_\star (\mathfrak{a}, L, M, \theta)$.

\medbreak

(ii) Take a sequence $(a^n)_n \subset {\cal C}$ which converges
towards $a \in {\cal C}$ with respect to the topology of $C^0_{\#}(\RR){}$.
For any $n$ let $\psi ^n,\ \rho ^n,\ \tilde\rho ^n,\ \lambda ^n,\
\tilde{a}^n = {\cal F}_{\mathfrak{a}}a^n$ constructed as in the
definition of the fixed point application. Similarly let $\psi ,\
\rho ,\ \tilde\rho ,\ \lambda ,\ \tilde{a} = {\cal F}_{\mathfrak{a}}a$.
The sequence $(\psi ^n )_n $ is bounded in $W^{1, \infty}(\RR)$
and therefore, by the Arzel\`a--Ascoli theorem we can extract a
subsequence converging in $C^0_{\#}(\RR){}$. Obviously the limit function
is $\psi$ and by the uniqueness of the limit we deduce that the
whole sequence $(\psi ^n)_n$ converges towards $\psi$ in $C^0_{\#}(\RR){}$.
In the same manner, since $\sup _n (\|\rho ^n \|_{L^\infty(\RR){}} +
\|\frac{\rmd}{\rmd x}{\rho ^n} \|_{L^\infty(\RR){}}) < +\infty $ we deduce that
$\rho ^n \rightarrow \rho$, $\tilde\rho ^n \rightarrow \tilde\rho $ in
$C^0_{\#}(\RR){}$. 

\smallbreak

The fact that $\lim _{n \rightarrow +\infty } \lambda ^n = \lambda$ stems from general spectrum continuity theorems~\cite{Kato66}, or can be directly deduced
from~\eref{Equ:lambda}.
Finally, as $\sup _n \|\tilde{a}^n\|_{W^{2,\infty}(\RR)} < +\infty $, 
we can extract
a subsequence $(\tilde{a}^{n_k})_k$ converging in~$C^1_{\#}(\RR)$ towards
some function~$\tilde{\tilde{a}}$. By passing to the limit with
respect to~$k$ in the weak formulation of $\tilde{a}^{n_k}$ we
obtain that the limit $\tilde{\tilde{a}}$ satisfies
$$
- \tilde{\tilde{a}}''(x) + \tilde\rho (x)\,
\tilde{\tilde{a}}(x) = \lambda\, \tilde{\tilde{a}}(x), \;\;x \in
(0,L),\;\;\tilde{\tilde{a}}(0) =
\tilde{\tilde{a}}(L),\;\tilde{\tilde{a}}'(0) =
\tilde{\tilde{a}}'(L).
$$
Moreover since $\tilde{a}^n \geq 0$, $\|\tilde{a}^n \|_{L^2(0,L){}} =
\mathfrak{a}$ for any $n$, we have $\tilde{\tilde{a}} \geq 0 $ and
$\|\tilde{\tilde{a}}\|_{L^2(0,L){}} = \mathfrak{a}$ and thus
$\tilde{\tilde{a}} = \tilde{a} = {\cal F}_{\mathfrak{a}} a$. By
the uniqueness of the limit we have $\lim _{n \rightarrow +\infty
} \tilde{a}^n = \tilde{a}$ in $C^1_{\#}(\RR)$.
\end{proof}

\medbreak

\noindent  We are now in position to prove our main result by
using the fixed point method.
\begin{theorem}
\label{MainResult} Assume that $\rho_\mathrm{b} \in L^1_{\#}(\RR)$, $\rho_\mathrm{b} \geq
0$, $\intx{\rho_\mathrm{b} (x)} = M$ and let $\theta $ be a positive real
number. For any $\mathfrak{a}>0$ there is at least one classical
solution $(\rho, a) \in C^1_{\#}(\RR) \times C^2_{\#}(\RR)$ for the
Boltzmann--Helmholtz equations~(\ref{Equ13}--\ref{Equ14})
satisfying $\rho \geq 0$, $\intx{\rho (x)} = M$, $a \geq 0$,
$\intx{a(x)^2} = \mathfrak{a}^2$.
\end{theorem}
\begin{proof}
Consider $\tilde{{\cal F}}_{\mathfrak{a}} = {\cal
F}_{\mathfrak{a}}|_{{\cal C}}$. The set ${\cal C}$ is convex and
compact in $C^0_{\#}(\RR){}$; by Proposition~\ref{FixedPointAppli} we
know that $\tilde{{\cal F}}_{\mathfrak{a}}({\cal C}) \subset {\cal
C}$ and that $\tilde{{\cal F}}_{\mathfrak{a}}$ is continuous with
respect to the topology of $C^0_{\#}(\RR)$. By the Schauder fixed point
theorem we deduce that there is a fixed point $a \in {\cal C}$. By
construction we have $a \geq 0$, $\intx{a(x)^2} = \mathfrak{a}^2$.
Consider now $\psi ,\ \rho,\ \tilde\rho,\ \lambda $ as in the definition
of $\tilde{{\cal F}}_{\mathfrak{a}} a$. Obviously $\lambda \geq 0,\ 
\rho \geq 0,\ \intx{\rho (x)} = M$ and we check easily that $(\rho
, a) \in C^1_{\#}(\RR) \times C^2_{\#}(\RR)$. Observe that $\lambda >0$. Indeed
we have
$$
\lambda = \frac{\intx{\{a'(x)^2 + \tilde\rho (x)\,a(x)^2\}}}{\intx{a(x)^2}}\geq \frac{\intx{\tilde\rho (x)\,a(x)^2}}{\intx{a(x)^2}}.
$$
If $\lambda = 0$ then $\tilde\rho (x)\,a(x)^2 = 0$ for any $x$, and since
by construction $\tilde\rho >0$ we deduce that $a = 0$ which
contradicts $\intx{a(x)^2} = \mathfrak{a}^2 > 0$. Consider now
$\omega = \sqrt{\lambda}>0$ and thus $(\rho ,a)$ is a solution 
of~(\ref{Equ13}--\ref{Equ14}).
\end{proof}

\section{Extensions.\label{exten}}
One could investigate the existence of ``non-linear
harmonics'' of the ``fundamental mode'' given by 
Theorem~\ref{MainResult}, \emph{i.e.},~solutions to~\eref{Equ14}
where $\omega^2$~is not the first eigenvalue of~$A_{\tilde\rho}$, but one of higher
rank. Unfortunately, it appears impossible to generalise the 
construction of~$\mathcal{F}_\mathfrak{a}$ to these eigenvalues. 
The reason is that, with periodic boundary conditions (unlike the 
Dirichlet, Neumann or Fourier b.c.), these eigenvalues may be double for 
some ``exceptional'' densities~$\tilde\rho$. For instance, if $\rho_\mathrm{b} 
=$~cst and $a =$~cst, then $\tilde\rho =$~cst and all eigenvalues except
the first one are double. There is apparently no way of defining a continuous 
mapping $\tilde\rho \mapsto \tilde{a}$ in the neighbourhood of the 
exceptional densities. Nevertheless, the existence of harmonics is 
very likely, as the eigenvalues are generically simple.

\medbreak

Another interesting extension is the case where the ion density is no longer given, but is also proportional to the Boltzmann factor. Let us denote by the subscript~1, resp.~2, the quantities relative to the electrons, resp.~ions; we introduce a new parameter~$\mu$ representing the electron/ion mass ratio. Then, we have
$f_1 \propto \rme^{-W_1/\theta_1}$ and $f_2 \propto \rme^{-W_2/\theta_2}$.
The energy~$W_1$ of one electron is given in section~\ref{deriv}; that of one ion is, according to the relativistic character:
\begin{eqnarray*}
W_2(x,p) = \case12\,\mu\, (p^2 + |a(x)|^2) - V(x),\;\;\mbox{(NR)},\\
W_2(x,p) = \mu^{-1}\,\sqrt{1+(\mu\,p)^2} + \case12\,\mu |a(x)|^2 - V(x),\;\;\mbox{(QR)},\\
W_2(x,p) = \mu^{-1}\,\sqrt{1+\mu^2\,(p^2 + |a(x)|^2)} - V(x),\;\;\mbox{(FR)}.
\end{eqnarray*}
We arrive at the following system:
\begin{eqnarray}
\label{Equ10:bip} V''(x) = \rho _2 (x) - \rho _1 (x),&& x \in \RR, \\
\label{Equ11:bip} - \omega ^2 a(x) - a''(x) = - (\tilde\rho_1(x) + \mu\, \tilde\rho_2(x))\, a(x),&\quad&  x \in \RR.
\end{eqnarray}
The arguments of sections~\ref{boltz} and~\ref{fixed} can be extended without bad surprises to this two-species model. However, the solutions corresponding to the first eigenvalue are not very interesting: one easily checks that all the functions $\rho_1,\ \rho_2,\ \tilde\rho_1,\ \tilde\rho_2,\ a$ are constant, and~$V\equiv0$.

\section{Concluding remarks.\label{concl}}
In this article, we have shown the existence of quasi-equilibrium solutions to the laser-plasma system, where the distribution function is Boltzmannian and the  electromagnetic variables are time-harmonic, at least at the fundamental frequency. The existence of solutions at higher frequencies 
is probable, both for one-species and two-species models. 
These solutions appear as generalisations of Vlasov--Poisson equilibria, but are clearly different from them as an electromagnetic wave is present. They can be viewed as a simple case of non-linear interaction between the electron plasma oscillations and the laser wave. The implicit relation (through the spectrum of the operator~$A_{\tilde\rho}$) between the frequency~$\omega$ and the space period~$L$ yields in the linear limit the dispersion relation for electromagnetic waves. 

\smallbreak

Quasi-equilibria can serve as references for analysing the dynamics of laser-plasma interaction, \emph{e.g.}~Raman and Brillouin scattering, which are among the most challenging issues to deal with in order to achieve controlled inertial confinement fusion.
Indeed, from a dynamical point of view, it should be noted that quasi-equilibria may be unstable, unlike the Vlasov--Poisson equilibria which are non-linearly stable, even under 1D Vlasov--Maxwell perturbations~\cite{CaLa06}.
These solutions may also serve as benchmarks for testing numerical codes, even though the numerical solution of the Boltzmann problem appears quite difficult when $u_\mathrm{b}$ and/or~$\psi$ feature large variations.

% ACKNOWLEDGEMENTS
\ack
The authors thank Pierre Bertrand for many stimulating discussions and for pointing out several useful references.

%%%%%%%%%%%%%
%           %
%  BIBLIO   %
%           %
%%%%%%%%%%%%%

\end{document}